\documentclass[11pt,reqno]{amsart}

\usepackage{epsfig,graphicx,color,mathrsfs}
\usepackage{graphicx}
\usepackage{amsmath,amssymb,amsthm,amsfonts}
\usepackage{amssymb}
\usepackage[english]{babel}
\usepackage{epsfig,graphicx,color,mathrsfs}
\usepackage[english]{babel}
\usepackage[left=2.7cm,right=2.7cm,top=3cm,bottom=3cm]{geometry}

\newtheorem{theorem}{Theorem}[section]
\newtheorem{proposition}[theorem]{Proposition}
\newtheorem{lemma}[theorem]{Lemma}

\newtheorem{remark}[theorem]{Remark}
\newtheorem{definition}[theorem]{Definition}

\newcommand{\R}{\mathbb R}

\newcommand{\e}{\varepsilon}

\newcommand{\io}{\int_{\Omega}}
\newcommand{\ho}{H^\circ}
\newcommand{\n}{\nabla}
\newcommand{\la}{\langle}
\newcommand{\ra}{\rangle}
\newcommand{\vf}{\varphi}
\newcommand{\Om}{\Omega}
\newcommand{\di}{\,\text{\rmfamily\upshape d}}
\newcommand{\br}{B_r(x_0)}
\newcommand{\brr}{B_{2r}(x_0)}
\newcommand{\bh}{\mathcal B^H}
\newcommand{\bho}{\mathcal B^{H^\circ}}
\newcommand{\ia}{\int_{A_R}}

\DeclareMathOperator{\Div}{div}

\title[Quasilinear anisotropic elliptic equations]{Hopf Lemma and regularity results for quasilinear anisotropic elliptic equations}

\author[D. Castorina]{Daniele Castorina$^*$}

\author[G. Riey]{Giuseppe Riey$^{**}$}

\author[B. Sciunzi]{Berardino Sciunzi$^{**}$}


\thanks{$^*$ Department of Mathematics, Computer Science and Natural Sciences,
John Cabot University, Via della Lungara 233, 00165 , Roma, Italy, E-mail: {\em dcastorina@johncabot.edu}}

\thanks{$^{**}$Dipartimento di Matematica e Informatica,
Universit\`a della Calabria,
Ponte Pietro Bucci 31B, I-87036 Arcavacata di Rende, Cosenza, Italy,
E-mail: {\em riey@mat.unical.it}, {\em sciunzi@mat.unical.it}}

\begin{document}

\begin{abstract}
We consider a class of quasi-linear anisotropic elliptic equations, possibly degenerate or singular, which are of interest in several applications such as computer vision and continuum mechanics.
We prove a Hopf Lemma as well as local and global regularity estimates for positive solutions, generalizing previous results known in the context of p-Laplacian equations.
\end{abstract}

\maketitle

\tableofcontents

\section{Introduction}\label{intro}

We consider quasilinear ellitic partial differential equations, in a possibly anisotropic medium. From the mathematical viewpoint, the anisotropy is responsible for a much richer geometric structure  than the usual Euclidean geometry. However, the main interest rests in concrete applications, since anisotropic media naturally arise in several real world phenomena.

In fact, anisotropic energies are widely used in computer vision (see for instance \cite{BaFe,EsOs,PeMa,We,YiZhSi}) and in continuum mechanics, in particular in presence of materials with distinct behavior with respect to different directions,
typically due to the crystalline microstructure of the medium
(see for instance \cite{BeNoPa,BeNoRi,BrRiSo,CaHo,Gu,Ta} and the references therein).

Our main results are a Hopf Lemma at the boundary, as well as local and global regularity estimates for positive solutions

Let us remark that, in a forthcoming paper, as a direct application of the results discussed above, we shall develop a moving plane procedure in the general context of this anisotropic and possibly singular/degenerate elliptic equations, in order to prove monotonicity and symmetry results in this Finsler geometry setting for positive solutions, both on bounded and unbounded domains, such as the whole space $\R^n$ or on half spaces.

For $n\geq 2$, let $\Omega \subset \R^n$ be a smooth bounded domain.
Let us consider the following \emph{Wulff type functional}:
\begin{equation}\label{wulff}
I(u)=\int_\Om\left[B(H(\n u))-F(u)\right]\di x\,,
\end{equation}
whose weak form of the Euler-Lagrange equation is given by
\begin{equation}\label{eq debole}
\int_\Om B'(H(\n u))\la\n H(\n u),\n \psi\ra =\int_\Om f(u)\psi, \, \forall \psi \in C^{1}_{c} (\Omega),
\end{equation}
(where $f=F'$)
as well as its strong form
\begin{equation}\label{eq forte}
-\Div\left(B'(H(\n u))\n H(\n u))\right)=f(u).
\end{equation}

We assume the following set of hypotheses on $B$, $H$ and $f$:
\begin{itemize}
\item[(i)] $B\in C^{3,\beta}_{loc}((0,+\infty))\cap C^1([0,+\infty))$, with $\beta\in (0,1)$
\item[(ii)] $B(0)=B'(0)=0,\quad B(t),B'(t),B''(t)>0\,\,\forall t\in(0,+\infty)$
\item[(iii)] there exist $p>1,k\in [0,1],\gamma>0,\Gamma>0$ such that:
\begin{equation}\label{propr B}
\begin{array}{c}
\gamma (k+t)^{p-2}t\leq B'(t)\leq \Gamma (k+t)^{p-2}t\\
\gamma (k+t)^{p-2}\leq B''(t)\leq \Gamma (k+t)^{p-2}
\end{array}
\end{equation}
for any $t>0$
\item[(iv)] $H\in C^{3,\beta}_{loc}(\R^n\backslash \{0\})$ is even and such that $H(\xi)>0\,\,\forall\xi\in\R^n\backslash\{0\}$
\item[(v)] $H(t\xi)=tH(\xi)\quad\forall\xi\in\R^n\backslash\{0\},\forall t>0$
\item[(vi)] $H$ is \emph{uniformly elliptic}, that means
the set $B_1^H:=\{\xi\in\R^n:\,H(\xi)<1\}$ is uniformly convex, i.e.
$$
\exists \lambda>0:\quad \la D^2 H(\xi)v,v\ra\geq\lambda |v|^2\quad\forall\xi\in\partial B_1^H,\,\forall v\in\n H(\xi)^\perp\,.
$$
\item[(vii)] $f$ is a positive continuous function on $[0,\infty)$ which is locally Lipschitz continuous on $(0,\infty)$
\item[(viii)] there exists $g\in C^0([0,+\infty))$ non-decreasing on $(0,\delta)$, $\delta>0$,
satisfying $g(0)=0$, $f+g\geq 0$ and either $g=0$ in $[0,d]$, $d>0$, or
$$
\int_0^\delta\frac{1}{L^{-1}(G(s))}ds=+\infty,
$$
where $G(s)=\displaystyle\int_0^\delta g(t)dt$ and $L(s)=sB'(s)-B(s)$.
\end{itemize}
The assumptions $(iv)-(v)-(vi)$ ensure that $H$ is a Finsler norm.
For $H(\xi)=|\xi|$ we get the usual Euclidean norm and, if we take $B(t)=\frac{t^p}{p}$,
the operator at left-hand side of \eqref{eq forte} becomes the usual $p$-Laplacian operator.

Let us observe that, under the above hypotheses, the natural space for the existence of solutions is
$W^{1,p} (\Om) \cap L^{\infty} (\Om)$.
However, better regularity holds in general.
If $k>0$ in hypothesis $(iii)$, then the operator is uniformly elliptic and, by standard elliptic regularity, solutions are classical and \eqref{eq forte} is satisfied.
On the other hand, for $k=0$, the operator becomes degenerate or singular, and solutions are not classical.
In fact, in \cite{CoFaVa} the authors show that to equation \eqref{eq forte}
it is possible to apply results in \cite{DB, To} to ensure that the solutions belong to
$C^{1,\alpha}(\Omega)\cap C^2(\Omega\backslash Z)$ for some $0 < \alpha < 1$,
where $Z$ denotes the \emph{critical set}, i.e. the set whet $\n u$ vanishes.
Moreover, if $\Omega$ is smooth, we can apply the regularity results in \cite{Lie}
to deduce that the solutions are in fact $C^1$ up to the boundary.
Thus, in order to cover the general case, throughout the paper we shall consider solutions belonging to $C^1 (\overline{\Om})$, which implies that the equation shall always be intended in the weak sense \eqref{eq debole}.
Anyhow, we will see that, as a consequence of our regularity results,
the critical set $Z$ is negligible and the strong equation \eqref{eq forte} shall be satisfied almost everywhere.

In \cite{DS1,DS2} were firstly introduced useful tools to get regularity and
qualitative properties (comparison principles, Harnack inequality, monotonicity and symmetries, etc.)
of equations involving the $p-$Laplace operator.
Then these techniques were widely developed to study several type of
more general equations and systems (for instance with lower order terms or singular data)
both in bounded and unbounded domains
(see for instance \cite{FaMoRiSc, LePoRi, MeRiSc, MoRiSc, RiSc, Sc1, Sc2} and the references quoted there).
Using this framework,  we prove local regularity estimates (Section~\ref{localreg})
for the solutions of our anisotropic elliptic quasilinear equations,
namely a weighted integral hessian estimate as well as the integrability of the inverse of the gradient.
For these kind of equations we also prove a Hopf type Lemma  (Section \ref{sect hopf}).
Thanks to this result, the local regularity estimates are then extended to the global case (Section~\ref{globalreg}).

\section{Notation and some geometrical tools}\label{notation}

For $a,b\in\R^n$ we denote by $a\otimes b$ the matrix whose entries
are $(a\otimes b)_{ij}=a_i b_j$. We remark that, for $v,w\in\R^n$,
there holds:
\begin{equation}\label{tensori}
\la a\otimes b\, v,w\ra=\la b,v\ra \la a,w\ra\,.
\end{equation}
Given an $n\times n$ matrix $A$, we set: $\displaystyle
|A|:=\sqrt{\sum_{i,j=1}^n |A_{ij}|^2}$.
For $x_0\in\R^n$ and $r>0$ we set $\br=\{x\in\R^n:|x-x_0|<r\}$.\\

We briefly recall some basic properties about Finsler geometry,
which is the main tool to study anisotropic problems.
We recall that Riemannian geometry is a particular case of the Finsler one
and, in fact, also in this more general framework it is possible to define
length of curves, geodesics, curvatures, etc.
For our purposes, we focus the attention on Finsler norms not depending on the position,
that means invariant by translations.

\begin{definition}\label{finsler norm}\rm
A function $H:\R^n\to [0,+\infty)$ is said a \emph{Finsler norm} if it is continuous,
even, convex and it satisfies:
\begin{equation}\label{finsler 1}
H(\lambda\xi)=|\lambda|H(\xi),\quad\forall\,\lambda\in\R,\quad\forall\,\xi\in\R^n
\end{equation}
and
\begin{equation}\label{finsler 2}
\exists\,c>0\,:\quad H(\xi)\geq c|\xi|\quad\forall\,\xi\in\R^n\,.
\end{equation}
\end{definition}

The \emph{dual norm} $\ho:\R^n\to [0,+\infty)$ is defined as:
$$
\ho(x)=\sup\{\la \xi,x\ra: H(\xi)\leq 1\}.
$$
It is easy to prove that $\ho$ is also a Finsler norm and it has the same regularity
properties of $H$. In particular $\ho$ satisfies \eqref{finsler 2} with $c^{-1}$ in place of $c$.
Moreover it follows that $(\ho)^\circ=H$.
For $r>0$ and $\overline x\in\R^n$ we define:
$$
\bh_r(\overline x)=\{x\in\R^n: H(x-\overline x)\leq r\}
$$
and
$$
\bho_r(\overline x)=\{x\in\R^n: \ho(x-\overline x)\leq r\}.
$$
For simplicity, when $\overline x=0$, we set: $\bh_r=\bh_r(0),\,\bho_r=\bho_r(0)$.
In literature $\bh_r$ and $\bho_r$ are also called ``Wulff  shape'' and
``Frank diagram'' respectively.

We remark that there holds:
\begin{equation}\label{propr finsler}
H(\n \ho(x))=1=\ho(\n H(x)).
\end{equation}
For more details on Finsler geometry see for instance \cite{BaChSh, BePa}.

\section{Local regularity estimates}\label{localreg}

The aim of this section is to present some integral regularity estimates for the hessian and for the inverse of the gradient of any (local) solution of \eqref{eq debole}. First of all, we shall need the following lemma about some structural bounds for the principal part of our divergence form operator.

\begin{lemma}\label{stime fondamentali}
There exist $\bar C_1,\bar C_2>0$ such that:
\begin{equation}\label{propr BH 1}
\la\left[
B''(H(\xi))\n H(\xi)\otimes\n H(\xi)+
B'(H(\xi))D^2 H(\xi)
\right]v, v\ra\geq
\bar C_1(k+|\xi|)^{p-2}|v|^2
\end{equation}
and
\begin{equation}\label{propr BH 2}
\left|
B''(H(\xi))\n H(\xi)\otimes \n H(\xi)+B'(H(\xi))D^2 H(\xi)
\right|\leq
\bar C_2(k+|\xi|)^{p-2}
\end{equation}
for any $\xi\in \R^n \backslash\{0\}$ and $v\in\R^n$.
\end{lemma}
\proof

See formulas (3.2) and (3.3) in \cite{CoFaVa}.
\endproof

Next, we shall be interested in the \textbf{linearized equation} of \eqref{eq debole} at any fixed solution $u$, which we can write as follows.
Set $Z=\{x\in\Om:\n u(x)=0\}$ and, for $\vf\in C^{\infty}_c(\Om\backslash Z)$, taking $\psi=\vf_i$ in \eqref{eq debole}, we get:
\begin{equation*}
\int_\Om B''(H(\n u))\la\n H(\n u),\n u_i\ra\la\n H(\n u),\n \vf\ra+ B'(H(\n u))\la D^2H(\n u)\n u_i,\n\vf\ra=
\int_\Om f'(u)\vf,
\end{equation*}
which, taking in account \eqref{tensori}, can also be written as:
\begin{equation}\label{linearizzato}
\int_\Om \la
\left[B''(H(\n u))\n H(\n u)\otimes \n H(\n u)+ B'(H(\n u))D^2H(\n u)\right]\n u_i,\n\vf\ra=
\int_\Om f'(u)\vf.
\end{equation}
Let us remark that we could make sense of \eqref{linearizzato} under several regularity hypotheses on the solution as well as on the test functions. However, even if we will not pursue such generality, let us point out that the right space for the linearization and a full spectral theory for this equation, in the singular/degenerate case when $k=0$, has been introduced in \cite{DS1},\cite{DS2} and completed in \cite{CES4}.

\begin{remark}
In the sequel, with a little abuse of notation, we will
denote by $\n u_i$ (and $u_{ij}$ respectively) the second derivatives
of $u$ outside $Z$ (thought extended equal to $0$ on $Z$).
Then at the end of the section we will recover the sufficient regularity
to ensure that actually these derivatives coincide with the distributional
second derivatives of $u$ in the whole of $\Omega$.
\end{remark}

We are now ready to prove one of our main regularity results, namely a local weighted integral estimate for the Hessian.

\begin{proposition}[Local Hessian estimate]\label{stima hessiano locale}\rm
Let $u\in C^1(\overline{\Omega})$ be a solution to \eqref{eq debole}.
Fix $x_0\in \Omega $ and $r>0$ such that $B_{2r}(x_0)\subset\Omega$. For $\beta\in [0,1)$ and $\gamma<n-2$ ($\gamma=0$ if $n=2$),
there holds:
\begin{equation}\label{eq stima hessiano locale componente}
\sup\limits_{y\in \Omega}\, \int_{\br}\frac{(k+|\n u|)^{p-2-\beta} | u_{ij}|  ^2}{|x-y|^\gamma}\di x \leq C
\end{equation}
and
\begin{equation}\label{eq stima hessiano locale}
\sup\limits_{y\in \Omega}\, \int_{\br}\frac{(k+|\n u|)^{p-2-\beta} | D^2 u |  ^2}{|x-y|^\gamma}\di x \leq C\,,
\end{equation}
where $C= C(x_0,r,\beta, \gamma,p,n,\|u\|_{W^{1,\infty}},f)$.
\end{proposition}
\proof
Let $G_\e:\R\to\R$ be defined as:
$$
G_\e(s)=\left\{
\begin{array}{ll}
s & \hbox{ if } |s| \geq 2 \e, \\
2\left[s- \e   \frac{s}{|s|}\right]  & \hbox{ if }  \e< |s|< 2\e, \\
0 & \hbox{ if } |s|\leq \e,
\end{array}
\right.
$$
and let $\psi$ be a cut-off
function such that
\begin{equation}\label{psi}
\psi\in C^\infty_c(B_{2r}(x_0))\quad \psi\equiv 1\ \mbox{ in } \ B_r(x_0)\qquad \mbox{and }\quad |\n \psi|\leq\frac{2}{r},
\end{equation}
with $2r <$ dist$(x_0,\partial \Omega)$.
Fix $\beta\in [0,1)$
and $\gamma<n-2$ (or $\gamma=0$ if $n=2$) and set:
\begin{equation}\label{test}
 \varphi (x)=T_\e(u_i(x)) K_\delta(|x-y|) \psi^2(x)
\quad  \mbox{ where } \quad T_\e(t)=\frac{G_\e(t)}{|t|^\beta}
\quad\hbox{and}\quad
K_\delta(t)=\frac{G_\delta(t)}{|t|^{\gamma+1}}\,.
\end{equation}

Substituting $\vf$ in \eqref{linearizzato}, we get:
\begin{eqnarray}
\nonumber && \int_\Om B''(H(\n u))\la\n H(\n u),\n u_i\ra^2 T'_\e(u_i) K_\delta(|x-y|)\psi^2\\
\nonumber &+&\int_\Om B''(H(\n u))\la\n H(\n u),\n u_i\ra \la\n H(\n u),\n K_\delta(|x-y|)\ra T_\e(u_i)\psi^2\\
\nonumber &+&\int_\Om B''(H(\n u))\la\n H(\n u),\n u_i\ra \la \n H(\n u),\n \psi\ra T_\e(u_i)K_\delta(|x-y|)2\psi\\
\nonumber &+&\int_\Om B'(H(\n u))\la D^2H(\n u)\n u_i,\n u_i\ra T'_\e(u_i)K_\delta(|x-y|)\psi^2\\
\nonumber &+&\int_\Om B'(H(\n u))\la D^2H(\n u)\n u_i,\n K_\delta(|x-y|)\ra T_\e(u_i)\psi^2\\
\nonumber &+&\int_\Om B'(H(\n u))\la D^2H(\n u)\n u_i,\n \psi\ra T_\e(u_i)K_\delta(|x-y|)2\psi\\
\label{stima hessiano 1} &=&\int_\Om f'(u) T_\e(u_i)K_\delta(|x-y|)\psi^2\,.
\end{eqnarray}

Recalling \eqref{tensori}, we have:
$$
\la \n H(\xi),\n v\ra^2=\la \n H(\xi)\otimes \n H(\xi)v,v\ra
\quad
\forall \xi\in\R^n,\,\forall v\in\R^n
$$
and
$$
\la \n H(\xi),v\ra
\la \n H(\xi),w\ra=
\la \n H(\xi)\otimes \n H(\xi) v,w\ra
\quad
\forall \xi\in\R^n,\,\forall v\in\R^n,\forall w\in\R^n.
$$
Hence by \eqref{propr BH 1} and \eqref{propr BH 2} we have:
\begin{equation}\label{stima hessiano 2}
B''(H(\n u))\la\n H(\n u),\n u_i\ra^2+
B'(H(\n u))\la D^2H(\n u)\n u_i,\n u_i\ra\geq
\bar C_1(k+|\n u|)^{p-2}|\n u_i|^2
\end{equation}
and
\begin{equation}\label{stima hessiano 2 bis}
\left| B''(H(\n u))\n H(\n u)\otimes \n H(\n u)+
B'(H(\n u)) D^2H(\n u)\right|\leq
\bar C_2(k+|\n u|)^{p-2}.
\end{equation}

By \eqref{stima hessiano 2} and Cauchy-Schwartz inequality we get:

\begin{eqnarray}
\nonumber && \bar C_1\int_\Om(k+|\n u|)^{p-2}|\n u_i|^2T'_\e(u_i) K_\delta(|x-y|)\psi^2\\
\nonumber &\leq& \int_\Om
\left|\la \left[B''(H(\n u))\n H(\n u)\otimes \n H(\n u)+B'(H(\n u)) D^2H(\n u)\right]\n u_i,\n K_\delta(|x-y|)\ra\right|
|T_\e(u_i)|\psi^2\\
\nonumber &+& \int_\Om
\left|\la \left[B''(H(\n u))\n H(\n u)\otimes \n H(\n u)+B'(H(\n u)) D^2H(\n u)\right]\n u_i,\n \psi\ra\right||T_\e(u_i)|K_\delta(|x-y|) 2\psi\\
\nonumber &+&\int_\Om |f'(u)||T_\e(u_i)|K_\delta(|x-y|)\psi^2\\
\nonumber &\leq& \int_\Om
\left|B''(H(\n u))\n H(\n u)\otimes \n H(\n u)+B'(H(\n u)) D^2H(\n u)\right||\n u_i||\n K_\delta(|x-y|)|
|T_\e(u_i)|\psi^2\\
\nonumber &+& \int_\Om
\left|B''(H(\n u))\n H(\n u)\otimes \n H(\n u)+B'(H(\n u)) D^2H(\n u)\right||\n u_i||\n \psi||T_\e(u_i)|K_\delta(|x-y|) 2\psi\\
\label{stima hessiano 3} &+&\int_\Om |f'(u)||T_\e(u_i)|K_\delta(|x-y|)\psi^2.
\end{eqnarray}

Taking into account \eqref{stima hessiano 2 bis}, by \eqref{stima hessiano 3} we infer:

\begin{eqnarray}\label{stima hessiano 4}
\nonumber && \bar C_1\int_\Om(k+|\n u|)^{p-2}|\n u_i|^2T'_\e(u_i) K_\delta(|x-y|)\psi^2\\
\nonumber &\leq& \bar C_2\int_\Om (k+|\n u|)^{p-2}|\n u_i||\n K_\delta(|x-y|)||T_\e(u_i)|\psi^2\\
\nonumber &+& \bar C_2\int_\Om (k+|\n u|)^{p-2}|\n u_i||\n \psi||T_\e(u_i)|K_\delta(|x-y|) 2\psi\\
&+&\int_\Om |f'(u)||T_\e(u_i)|K_\delta(|x-y|)\psi^2.
\end{eqnarray}

Since $\Om$ is bounded, we have that $\int_\Omega |x-y|^{-s} \di x$ is uniformly bounded for any $s<n$. In particular, since $\gamma<n-2$, this is true both for $s =\gamma$ and $s=\gamma+1$.
Therefore, for fixed $\e>0$, we can send $\delta$ to $0$ in \eqref{stima hessiano 4} and, recalling the definition of $K_\delta$,
by dominated convergence we get:

\begin{eqnarray}\label{stima hessiano 5}
\nonumber && \int_\Om \frac{(k+|\n u|)^{p-2}|\n u_i|^2T'_\e(u_i) \psi^2}{|x-y|^\gamma}\\
\nonumber &\leq& C\int_\Om \frac{(k+|\n u|)^{p-2}|\n u_i||T_\e(u_i)|\psi^2}{|x-y|^{\gamma+1}}\\
\nonumber &+& C \int_\Om \frac{(k+|\n u|)^{p-2}|\n u_i||\n \psi||T_\e(u_i) |\psi}{|x-y|^\gamma}\\
&+&C \int_\Om \frac{|f'(u)||T_\e(u_i)| \psi^2}{|x-y|^\gamma}.
\end{eqnarray}

We remark that, if $n=2$, the first term in the sum at the right hand-side of \eqref{stima hessiano 5} is equal to $0$
because $\n K_\delta=0$ if $\gamma=0$. If instead $n\geq 3$, recalling that $G_\e$ is an odd function and that $\gamma<n-2$,
for a $0<\theta<1$ we have:
\begin{equation}\label{stima I1}
\int_\Om \frac{(k+|\n u|)^{p-2}|\n u_i||T_\e(u_i)|\psi^2}{|x-y|^{\gamma+1}}
\leq
\theta \int_\Omega\frac{(k+|\n u|)^{p-2}|\n u_i|^2}{|x-y|^\gamma |u_i|^\beta}
\frac{G_\e(u_i)}{u_i} \psi^2+ C.
\end{equation}

Since $|\n\psi|\leq\frac{2}{r}$, we have:
\begin{equation}\label{stima I2}
\int_\Om \frac{(k+|\n u|)^{p-2}|\n u_i||\n \psi||T_\e(u_i) |\psi}{|x-y|^\gamma}
\leq \theta \int_\Omega\frac{(k+|\n u|)^{p-2}|\n u_i|^2 }{|x-y|^\gamma |u_i|^\beta }
\frac{G_\e(u_i)}{u_i}\psi^2+C.
\end{equation}

Recalling that $\beta\in[0,1)$ and that $|T_\e(u_i)|\leq |u_i|^{1-\beta}$,
since $f$ is Lipschitz and $u$ is bounded we get:
\begin{equation}\label{stima I3}
\int_\Om \frac{|f'(u)||T_\e(u_i) | \psi^2}{|x-y|^\gamma}
\leq C\int_{\Omega} \frac{1}{|x-y|^\gamma}\di x\leq C.
\end{equation}

For $s>0$ we have
$$
T'_\e(s)=\frac{1}{|s|^\beta}\left[G'_\e(s)-\beta\frac{G_\e(s)}{s}\right]
$$
and by \eqref{stima hessiano 5}, \eqref{stima I1}, \eqref{stima I2} and \eqref{stima I3} we get:
\begin{equation}\label{eq stima hess 6}
\int_\Omega \frac{(k+|\n u|)^{p-2}|\n u_i|^2}{|u_i|^\beta |x-y|^\gamma}
\left(
G'_\e(u_i)-(\beta+\vartheta)\frac{G_\e(u_i)}{u_i}\right)\psi^2\leq C.
\end{equation}
We now choose $\theta$ small enough so that $\beta+\theta<1$,
so that $G'_\e(u_i)-(\beta+\vartheta)\frac{G_\e(u_i)}{u_i}$ is positive.
By definition of $G_\e$ it follows that, for any $s >0$,
$G'_\e(s)-(\beta+\vartheta)\frac{G_\e(s)}{s}$ tends to $1-(\beta+\theta)$ as
$\e$ goes to $0$, and hence by Fatou's Lemma we get:
\begin{equation}\label{eq stima hess 7}
\int_{\Omega\setminus\{u_i=0\}} \frac{(k+|\n u|)^{p-2}|\n u_i|^2}{|u_i|^\beta |x-y|^\gamma}\psi^2\leq C
\end{equation}
and, since $(k+|\n u|)^\beta\geq |\n u|^\beta\geq |u_i|^\beta$, we have:
\begin{equation}\label{eq stima hess 8}
\int_{\Omega\setminus\{u_i=0\}} \frac{(k+|\n u|)^{p-2-\beta}|\n u_i|^2}{|x-y|^\gamma}\psi^2\leq C
\end{equation}
whence, recalling that $u_{ij} = 0$ on $\{u_i=0\}\cap\left(\Omega\setminus Z\right)$, we see that:
\begin{equation}\label{eq stima hess 9}
\int_{\Omega\setminus Z} \frac{(k+|\n u|)^{p-2-\beta}|\n u_i|^2}{|x-y|^\gamma}\psi^2\leq C\,,
\end{equation}
where $C$ depends on $x_0,r, n,p,\beta,\gamma, f, H, \|u\|_{W^{1,\infty}}$
but it does not depend on $y$.
By the properties of $\psi$ it follows:
$$
\sup\limits_{y\in \Omega}\,\, \int_{B_r(x_0)\setminus Z}\frac{(k+|\n u|)^{p-2-\beta} |D^2 u|^2}{|x-y|^\gamma}\leq C
$$
and by standard arguments (see for instance \cite{DS1},\cite{Sc1},\cite{Sc2}) we get the thesis.

 In fact, inspired by Stampacchia's Theorem, we first
 extend the (generalized) second derivatives of $u$ to be zero over the critical set $Z$ and we can actually write that
\begin{equation}\label{stima hess piena}
\sup\limits_{y\in \Omega}\,\, \int_{B_r(x_0)}\frac{(k+|\n u|)^{p-2-\beta} |D^2 u|^2}{|x-y|^\gamma}\leq C
\end{equation}
Such an estimate then allows us to conclude that the extended generalized derivatives are the effective distributional derivatives (see e.g. \cite{DS1} for details).

\endproof

At this point, before proving the integrability of the inverse of the gradient, we shall need another structural estimate.

\begin{lemma}\label{stima BH}
There exists $C>0$ such that:
\begin{equation}\label{eq stima BH}
|B'(H(\xi))|\leq C \left(k+|\xi|\right)^{p-1}.
\end{equation}
\end{lemma}
\proof
Since $H$ is a norm equivalent to the euclidean one, there exist $\lambda_1,\lambda_2>0$ such that:
\begin{equation}\label{norma equivalente}
\lambda_1|\xi|\leq H(\xi)\leq \lambda_2 |\xi|\quad\quad\forall \xi\in\R^n.
\end{equation}
Hence by \eqref{propr B} we infer:
\begin{equation}\label{stima BH 1}
B'(H(\xi))\leq C_2(k+H(\xi))^{p-2}H(\xi)\leq C_2(k+H(\xi))^{p-2}\lambda_2|\xi|.
\end{equation}
Moreover we have:
\begin{equation}\label{stima BH 2}
\left\{
\begin{array}{ll}
(k+H(\xi))^{p-2}\leq (k+\lambda_2|\xi|)^{p-2} & \hbox{ if } p\geq 2\\
(k+H(\xi))^{p-2}\leq (k+\lambda_1|\xi|)^{p-2} & \hbox{ if } p<2.
\end{array}
\right.
\end{equation}
Let us consider first the case $p\geq 2$.
For $t>0$, if $\lambda_2\leq 1$, we have:
\begin{equation}\label{stima BH 3}
k+\lambda_2 t\leq k+t,
\end{equation}
while, if $\lambda_2>1$, we have:
\begin{equation}\label{stima BH 4}
k+\lambda_2 t=\lambda_2\left(\frac{k}{\lambda_2}+t\right)\leq \lambda_2(k+t).
\end{equation}
By \eqref{stima BH 1}, the first equation in \eqref{stima BH 2}, \eqref{stima BH 3}
and \eqref{stima BH 4}, we get:
\begin{equation}\label{stima BH 5}
B'(H(\xi))\leq C_2(k+\lambda_2|\xi|)^{p-2}\lambda_2|\xi|\leq C_2(k+\lambda_2|\xi|)^{p-1}
\leq C_2\max\{1,\lambda_2\}^{p-1}(k+|\xi|)^{p-1}.
\end{equation}
Let us now consider the case $p<2$.
For $t>0$, if $\lambda_1\geq 1$, we have:
\begin{equation}\label{stima BH 6}
k+\lambda_1 t\geq k+t,
\end{equation}
while, if $\lambda_1<1$, we have:
\begin{equation}\label{stima BH 7}
k+\lambda_1 t=\lambda_1\left(\frac{k}{\lambda_1}+t\right)> \lambda_1(k+t).
\end{equation}
Hence by \eqref{stima BH 1}, the second equation in \eqref{stima BH 2}, \eqref{stima BH 6}
and \eqref{stima BH 7}, we get:
\begin{equation}\label{stima BH 8}
\begin{split}
B'(H(\xi)) &\leq \frac{C_2\lambda_2}{\lambda_1}(k+\lambda_1|\xi|)^{p-2}\lambda_1|\xi|\\
&\leq \frac{C_2\lambda_2}{\lambda_1}(k+\lambda_1|\xi|)^{p-1}\leq
\frac{C_2\lambda_2}{\lambda_1}\min\{1,\lambda_1\}^{p-1}(k+|\xi|)^{p-1}.
\end{split}
\end{equation}
\endproof

We are now in position to state and prove our second main regularity result, which deals with the local integrability of the inverse of the weight.

\begin{proposition}[Local estimate of the weight]\label{inverso peso locale}\rm
Let $u\in C^{1}(\overline{\Om})$ be a solution to \eqref{eq debole}.
Fix $t\in [0,p-1)$ and $\gamma<n-2$ ($\gamma=0$ if $n=2$).
Then, for any $\Omega'\subset\subset\Omega $ there exists $C$ such that
\begin{equation}\label{stima peso locale}
\sup\limits_{y\in \Omega}\,\, \int_{\Omega'}\frac{1}{(k+|\n  u|)^t |x-y|^\gamma} \di x \leq C\,,
\end{equation}
where $C=C(\Omega',t,\gamma,n,p,\|u\|_{W^{1,\infty}}, f)$.
\end{proposition}
\proof
We first prove inequality \eqref{stima peso locale}
on balls, then the thesis will follow by a covering argument.
For $x_0\in \Omega $ we choose $r>0$ such that $B_{2r}(x_0)$
is contained in $\Omega $.
Let $\psi$ and $K_\delta$ be defined as in Proposition \ref{stima hessiano locale}.
For $t=p-2+\beta<p-1$ and $\e > 0$, we consider in \eqref{eq debole} the test function:
$$
\vf=\frac{1}{(k+|\n u|)^t+\e}K_\delta(|x-y|)\psi^2
$$
and, noticing that $f(u(x)) \geq C(x_0)>0$ for any $x \in \brr$, we
get:
\begin{eqnarray}\label{peso 1}
&& C(x_0)\int_{\brr}\frac{1}{(k+|\n u|)^t+\e}K_\delta(|x-y|)\psi^2\\
\nonumber &\leq&
   -t\int_{\brr}\frac{B'(H(\n u))(k+|\n u|)^{t-1}}{[(k+|\n u|)^t+\e]^2}
   \left\la\n H(\n u),\frac{\n u}{|\n u|}D^2 u\right\ra K_\delta(|x-y|)\psi^2 \\
\nonumber &+& \int_{\brr}\frac{B'(H(\n u)) \psi^2}{(k+|\n u|)^t+\e}
   \la \n H(\n u),\n K_\delta(|x-y|)\ra\\
\nonumber &+& \int_{\brr}\frac{B'(H(\n u)) 2\psi K_\delta(|x-y|)}{(k+|\n u|)^t+\e}
   \la \n H(\n u),\n \psi\ra.
\end{eqnarray}
Recalling that $H$ is $1$-homogeneous, we have that $\n H$ is
$0$-homogeneous and hence we have:
$$
\n H(\xi)=\n H\left(|\xi|\frac{\xi}{|\xi|}\right)=\n H\left(\frac{\xi}{|\xi|}\right)\quad\forall \xi\in\R^n.
$$
Since $\n H$ is continuous, we infer that there exists $M>0$ such that:
\begin{equation}\label{bound grad H}
|\n H(\xi)|\leq M\quad\forall\xi\in\R^n.
\end{equation}
By \eqref{peso 1} and \eqref{bound grad H} and Lemma \ref{stima BH} we argue:
\begin{eqnarray}\label{peso 2}
&& \int_{\brr}\frac{1}{(k+|\n u|)^t+\e}K_\delta(|x-y|)\psi^2\\
\nonumber &\leq&
   C\int_{\brr}\frac{(k+|\n u|)^{p-1}\cdot (k+|\n u|)^{t-1}}{[(k+|\n u|)^t+\e]^2}
   |D^2 u||K_\delta(|x-y|)|\psi^2 \\
\nonumber &+& C\int_{\brr}\frac{(k+|\n u|)^{p-1}}{(k+|\n u|)^t+\e}
   |\n K_\delta(|x-y|)|\psi^2\\
\nonumber &+& C\int_{\brr}\frac{(k+|\n u|)^{p-1}}{(k+|\n u|)^t+\e}
   |K_\delta(|x-y|)||\n \psi|\psi.
\end{eqnarray}
Sending $\delta$ to $0$, by dominated convergence we get:
\begin{eqnarray}\label{peso 3}
&& \int_{\brr}\frac{1}{(k+|\n u|)^t+\e}\frac{1}{|x-y|^\gamma}\psi^2\\
\nonumber &\leq&
   C\int_{\brr}\frac{(k+|\n u|)^{p-1}\cdot (k+|\n u|)^{t-1}}{[(k+|\n u|)^t+\e]^2}
   |D^2 u|\frac{1}{|x-y|^\gamma}\psi^2 \\
\nonumber &+& C\int_{\brr}\frac{(k+|\n u|)^{p-1}}{(k+|\n u|)^t+\e}
   \frac{1}{|x-y|^{\gamma+1}}\psi^2\\
\nonumber &+& C\int_{\brr}\frac{(k+|\n u|)^{p-1}}{(k+|\n u|)^t+\e}
   \frac{1}{|x-y|^\gamma}|\n \psi|\psi.
\end{eqnarray}
Recalling that $t=p-2+\beta$, using Proposition \ref{stima hessiano locale},
for $0 <\theta <1$ we get:
\begin{eqnarray}\label{stima pezzo 1}
&& \int_{\brr}\frac{(k+|\n u|)^{p-1}\cdot (k+|\n u|)^{t-1}}{[(k+|\n u|)^t+\e]^2}
   |D^2 u|\frac{1}{|x-y|^\gamma}\psi^2 \\
\nonumber &\leq& \theta\int_{\brr}
  \frac{(k+|\n u|)^{3t}}{[(k+|\n u|)^t+\e]^4 |x-y|^\gamma}\psi^2
  +
  \frac{1}{4\theta}\int_{\brr}
  \frac{(k+|\n u|)^{p-2-\beta}| D^2 u |^2}{|x-y|^\gamma}
  \psi^2\\
\nonumber &\leq& \theta\int_{\brr}
  \frac{1}{[(k+|\n u|)^t+\e]|x-y|^\gamma}\psi^2
  +\frac C\theta.
\end{eqnarray}
Since $t<p-1$, $\frac{(k+|\n u|)^{p-1}}{(k+|\n u|)^t+\e}$ is bounded
and hence we have:
\begin{equation}\label{stima pezzo 2}
\int_{\brr}\frac{(k+|\n u|)^{p-1}}{(k+|\n u|)^t+\e}
   \frac{1}{|x-y|^{\gamma+1}}\psi^2\leq
   C\int_{B_{2r}(x_0)}\frac{1}{|x-y|^{\gamma+1}}\psi^2\leq C.
\end{equation}
Moreover by properties of $\psi$ it follows:
\begin{equation}\label{stima pezzo 3}
\int_{\brr}\frac{(k+|\n u|)^{p-1}}{(k+|\n u|)^t+\e}
   \frac{1}{|x-y|^\gamma}|\n \psi|\psi
\leq
\frac{C}{r}\int_{B_{2r}(x_0)}\frac{1}{|x-y|^\gamma}\psi^2\leq C.
\end{equation}
Choosing $\theta$ small enough, by \eqref{peso 3}, \eqref{stima pezzo 1},
\eqref{stima pezzo 2} and \eqref{stima pezzo 3} we get:
\begin{equation}\label{peso 4}
 \int_{\brr} \frac{1}{(k+|\n u|)^t+\e} \frac{1}{|x-y|^\gamma}\psi^2\leq C
\end{equation}
and, sending $\e$ to $0$, by Fatou's Lemma we get the thesis.
\endproof

\section{Hopf Lemma}\label{sect hopf}

In this section we shall prove a Hopf Lemma for solutions
of our anisotropic quasilinear elliptic equation.
We begin with a couple of structural bounds from below for our principal part
and then we state a weak comparison principle for a solution and a super-solution of an equation related to ours.

\begin{lemma}\label{lemma disug basso}
For $x\in\R^n\backslash\{0\},\,y\in\R^n$ and $p>1$ there exists $C>0$ such that:
\begin{equation}\label{disug basso 1}
\sum_{i,j=1}^{n}\frac{\partial}{\partial
x_i}\left(B'(H(x))H_j(x)\right)y_iy_j\geq C|x|^{p-2}|y|^2.
\end{equation}
and
\begin{equation}\label{disug basso 2}
\la B'(H(x))\n H(x)-B'(H(y))\n H(y),x-y\ra\geq
C\left(|x|+|y|\right)^{p-2}|x-y|^2.
\end{equation}
\end{lemma}
\proof First we remark that
$$
\frac{\partial}{\partial x_i}\left(B'(H(x)H_j(x)\right)=
B''(H(x))H_i(x)H_j(x)+B'(H(x))H_{ji}(x)
$$
and hence \eqref{disug basso 1} immediately follows by \eqref{propr BH 1}. Moreover, assuming $|y|\geq |x|$, using \eqref{disug basso 1}, we
have:
\begin{equation}\label{disug basso 2a}
B'(H(x))H_j(x)-B'(H(y))H_j(y)=
\int_0^1\sum_{i=1}^n\frac{\partial}{\partial
x_i}\left.\left(B'(H(x))H_j(x)\right)\right|_{y+t(x-y)}(x_i-y_i)dt
\end{equation}
and hence
\begin{eqnarray}\label{disug basso 2b}
&&\la B'(H(x)\n H(x)-B'(H(y))\n H(y),x-y\ra=\\
\nonumber &=&\int_0^1
  \sum_{i=1}^n
  \frac{\partial}{\partial x_i}\left.\left(B'(H(x))H_j(x)\right)\right|_{y+t(x-y)}(x_i-y_i)(x_j-y_j)dt\geq\\
\nonumber &\geq&C|x-y|^2\int_0^1|y+t(x-y)|^{p-2}dt.
\end{eqnarray}
For $t\in [0,1]$ it holds $|y+t(x-y)|\leq |x|+|y|$ and hence, if
$1<p<2$, then \eqref{disug basso 2} immediately follows. Else, if $p>2$, then we need to prove that
\begin{equation}\label{disug basso integrale}
\int_0^1|y+t(x-y)|^{p-2}dt\geq C (|x|+|y|)^{p-2}.
\end{equation}
We recall that we are assuming $|y|\geq |x|$ and hence, if $|x-y|\leq \frac{|y|}{2}$, then for $0<t<1$ it holds
$$
|y+t(x-y)|\geq|y|-|x-y|\geq\frac{|y|}{2}\geq \frac{|x|+|y|}{4},
$$
from which \eqref{disug basso integrale} follows. If instead
$|x-y|>\frac{|y|}{2}>0$, we set $t_0=\frac{|y|}{|x-y|}$ so that
$t_0\in (0,2)$ and we have:
\begin{eqnarray}\label{disug basso intermedia}
|y+t(x-y)| &\geq& ||y|-t|x-y||=|t-t_0||x-y|\geq\\
\nonumber &\geq& |t_0-t|\frac{|y|}{2}\geq |t_0-t|\frac{|x|+|y|}{4}
\end{eqnarray}
Recalling that $p>2$, we have $\int_0^1|t_0-t|^{p-2}dt\geq C$ and
hence \eqref{disug basso intermedia} implies \eqref{disug basso
integrale}.
\endproof

Thanks to the above lemma we can now deal with the following weak comparison principle. 

\begin{lemma}\label{lemma principio massimo}
Let $u,v\in C^1(\overline\Om)$ satisfy:
\begin{equation}\label{eq per princ max}
\left\{
\begin{array}{ll}
-\hbox{div}\left(B'(H(\n u))\n H(\n u)\right)+g(u)\geq 0 & \hbox{ in } \Om\\
 \\
-\hbox{div}\left(B'(H(\n v))\n H(\n v)\right)+g(v)=0 & \hbox{ in } \Om\\
 \\
v\leq u & \hbox{ on } \partial\Om\,,
\end{array}
\right.
\end{equation}
where $g$ satisfies the assumption $(viii)$ given in Section \ref{intro}.
Then there holds:
$$
v\leq u \hbox{ in }\Omega.
$$
\end{lemma}
\proof By the weak formulation of \eqref{eq per princ max} we get:
\begin{equation}\label{diseq princ max}
\io \la B'(H(\n v))\n H(\n v)-B'(H(\n u))\n H(\n u) ,\n\psi\ra \leq
\io (-g(v)+g(u))\psi.
\end{equation}
Taking $\psi=(v-u)^+$ as test function in \eqref{diseq princ max},
using \eqref{disug basso 2} and, recalling that $-g$ is
non-increasing, we infer:
\begin{equation}\label{diseq princ max intermedia}
\io \left(|\n u|+|\n v|\right)^{p-2}|\n (v-u)^+|^2\leq
\io(-g(v)+g(u))(v-u)^+\leq 0.
\end{equation}
If $p>2$, it follows that $\n u=\n v$ a.e. in $\Om$ and hence
$(v-u)^+$ is constant and, since $(v-u)^+=0$ on the boundary, we infer
that $(v-u)^+=0$. If $1<p<2$, since $u,v\in C^1(\overline\Om)$
\eqref{diseq princ max intermedia} implies:
$$
\io|\n (v-u)^+|^2\leq 0
$$
and hence as above we get again the thesis.
\endproof

\begin{remark}\label{remark principio massimo}
If $u$ satisfies \eqref{eq forte}, by property $(viii)$ of $f$ we infer that
$u$ satisfies also the first inequality in \eqref{eq per princ max}.
\end{remark}

At this point, in order to prove a Hopf Lemma at the boundary for any positive solution,
exactly as in the classical semilinear case, we shall need to construct a radial barrier from below defined in an annulus.

\begin{lemma}\label{esistenza radiale}
For $R>0$ and $\overline x\in\R^n$ we consider the annulus
$$
A_R(\overline x)=\bho_R(\overline
x)\backslash\overline{\bho_{\frac{R}{2}}(\overline x)}.
$$
Let $g$ satisfy the assumption $(viii)$ in Section \ref{intro}. Then
for every $m>0$ there exists a non-negative function $v\in
C^1(\overline {A_R})$ satisfying:
\begin{equation}\label{propr v}
\left\{
\begin{array}{ll}
-\Div\left(B'(H(\n v))\n H(\n v))\right)+g(v)=0 & \hbox{ in } A_R\\
v=0 & \hbox{ on } \partial\bho_R\\
v=m & \hbox{ on } \partial\bho_{\frac{R}{2}}\\
\frac{\partial v}{\partial \nu}>0 & \hbox{ on } \partial\bho_R\,,
\end{array}
\right.
\end{equation}
where $\nu$ denotes the inner unit normal to $\bho_R$.
\end{lemma}
\proof We look for radial solutions. The word ``\emph{radial}''
has to be understood in the Finsler framework: from now on we will say  that $v$ is radial if there exists $\overline x$ such
that $v$ is constant on the Finsler balls $\bho_R(\overline x)$ for any $R>0$. In
the sequel for simplicity we assume $\overline x=0$ and we set
$A_R=A_R(0)$. If we assume that $v$ is radial (with $\overline
x=0$), then there exists $w:[0,+\infty)\to\R$ such that:
\begin{equation}\label{def vw}
v(x)=w(\ho(x)).
\end{equation}
Using \eqref{def vw} in \eqref{eq debole} (where we take test
functions of the form $\psi(x)=\vf(\ho(x))$ with $\vf:\R\to\R$), we
get:
\begin{eqnarray}\label{eq debole radiale 1}
\nonumber &&\ia B'\left(H(w'(\ho(x))\n\ho(x))\right)
\left\la \n H\left(w'(\ho(x))\n\ho(x)\right),\vf'(\ho(x))\n\ho(x)\right\ra\di x+\\
&&+\ia g(w(\ho(x)))\vf(\ho(x))\di x=0.
\end{eqnarray}
Using the positive $1$-homogeneity of $H$ (and hence the
$0$-homogeneity of $\n H$), we infer:
\begin{eqnarray}\label{eq debole radiale 2}
\nonumber &&\ia B'\left(|w'(\ho(x))|H(\n\ho(x))\right)\hbox{sign}\left(w'(\ho(x))\right) \vf'(\ho(x))
\left\la \n H\left(\n\ho(x)\right),\n\ho(x)\right\ra\di x+\\
&&+\ia g(w(\ho(x)))\vf(\ho(x))\di x=0.
\end{eqnarray}
Recalling that $\la\n H(\xi),\xi\ra=H(\xi)$ and using \eqref{propr finsler}, \eqref{eq debole radiale 2} becomes:
\begin{equation}\label{eq debole radiale 3}
\ia B'\left(|w'(\ho(x))|\right)\hbox{sign}\left(w'(\ho(x))\right)\vf'(\ho(x))\di x+ \ia
g(w(\ho(x)))\vf(\ho(x))\di x=0.
\end{equation}
We now consider on $A_R$ a system of Finsler polar coordinates in
the following sense. Since $\bho_1$ is convex, we can find a map
$s:U\subset\R^{n-1}\to\R^n$ such that $\overline{s(U)}=\bho_1$.
Therefore for $\theta=(\theta_1,\cdots,\theta_{n-1})\in U$ we have
$\ho(s(\theta))=1$. For $\rho>0$ we have $\ho(\rho s(\theta))=\rho$
and hence $\rho s(\theta)$ belongs to $\bho_\rho$. This allows us to
consider the transformation $S:\R^n\to\R^n$ defined as:
$$
S(\rho,\theta)=\rho s(\theta)=(\rho s^1(\theta),\cdots,\rho s^n(\theta))
$$
and the change of variables:
$$
x=S(\rho,\theta)
$$
with $\rho\in \left[\frac{R}{2},R\right]$ and $\theta\in U$. Denoted
by $DS$ the Jacobian matrix of $S$, we set
$J(\rho,\theta)=|\hbox{det}DS(\rho,\theta)|$ and we have that there
exists a suitable function $\Gamma(\theta)>0$ such that:
\begin{equation}\label{jacobiano}
J(\rho,\theta)=\rho^{n-1}\Gamma(\theta).
\end{equation}
Using this change of variables, \eqref{eq debole radiale 3} becomes:
\begin{eqnarray}\label{eq debole radiale 4}
&&\int_U\Gamma(\theta)\di \theta_1\cdots\di \theta_{n-1}
 \int_{\frac{R}{2}}^R B'\left(|w'(\rho)|\right)\hbox{sign}\left(w'(\rho)\right)\vf'(\rho)\rho^{n-1}\di \rho+\\
\nonumber &+&
 \int_U\Gamma(\theta)\di \theta_1\cdots\di \theta_{n-1}
 \int_{\frac{R}{2}}^R g(w(\rho))\vf(\rho)\rho^{n-1}\di \rho=0
\end{eqnarray}
and hence:
\begin{equation}\label{eq debole radiale 5}
\int_{\frac{R}{2}}^R
B'\left(|w'(\rho)|\right)\hbox{sign}\left(w'(\rho)\right)\vf'(\rho)\rho^{n-1}\di \rho+
\int_{\frac{R}{2}}^R g(w(\rho))\vf(\rho)\rho^{n-1}\di \rho=0,
\end{equation}
whose form is:
\begin{equation}\label{eq debole radiale 6}
-\left(B'\left(|w'(\rho)|\right)\hbox{sign}\left(w'(\rho)\right)\rho^{n-1}\right)'+
g(w(\rho))\rho^{n-1}=0.
\end{equation}
If we are interested in Finsler radial solutions in the annulus
$A_R(\overline x)$, we can choose $\rho=\ho(x-\overline x)$ and for
$\rho\in \left[0,\frac{R}{2}\right]$ we consider
$q(\rho):=(R-\rho)^{n-1}$. In such a way $q$ satisfies the
assumptions required to apply Proposition 4.2.1 and Proposition 4.2.2 in \cite{PuSe},
which states that the problem:
\begin{equation}\label{eq debole radiale 7}
\left\{
\begin{array}{ll}
\left(\Phi\left(w'(\rho)\right)q(\rho)\right)' - g(w(\rho))q(\rho)=0 & \hbox{ in } (0,T)\\
w(0)=0,\quad w(T)=m>0,\quad w'(0)>0
\end{array}
\right.
\end{equation}
with $\Phi(t) = B'(t)\hbox{sign}(t)$, $q$ as defined above, $f=g$ and $T=R/2$ in our case, admits a $C^1$ solution satisfying $w'>0$.
\endproof

We close this section with the following theorem, which states that
the Hopf Lemma holds also for our equation. Let us point out that this can be seen an extension to the anisotropic setting of a classical Hopf Lemma for quasilinear equations by Vazquez \cite{Vaz}. 

\begin{theorem}[Hopf Lemma]\label{hopf}
Let $u$ be a $C^1(\overline\Om)$ solution of \eqref{eq debole} satisfying $u>0$ in
$\Om$ and $u=0$ on $\partial\Om$. Let $y\in\partial\Om$ such that
$\Om$ satisfies the interior sphere condition\footnote{ We remark
that the ``\emph{interior sphere condition}'' in the Euclidean sense
is equivalent to that in the Finsler geometry framework, that means,
if $\Om$ satisfies this condition with classical euclidean balls,
then it satisfies the same condition also if we consider the Finsler
balls $\bh$ or $\bho$. } at $y$ and let $\nu$ denote the inner unit
normal to $\partial\Om$ at $y$. Then there holds:
$$
\frac{\partial u}{\partial \nu}(y)>0.
$$
\end{theorem}
\proof Let $R>0$ be small enough such that there exists $\overline
x\in\Om$ such that $\bho_R(\overline x)\subset\Om$ and
$y\in\partial\bho_R(\overline x)$. Let $v$ be the function given by
Lemma \eqref{esistenza radiale}, which clearly satisfies the second
equality in \eqref{eq per princ max}.
Then by Lemma \ref{lemma principio massimo} and Remark \ref{remark principio massimo}
we have $u\geq v$ in $A_R(\overline x)$ and hence
$\displaystyle\frac{\partial u}{\partial
\nu}(y)\geq\displaystyle\frac{\partial v}{\partial \nu}(y)>0$.
\endproof

\section{Global regularity estimates}\label{globalreg}

In this section, thanks to the Hopf Lemma we just proved, we will see how it is possible to easily extend the local Hessian regularity result of Section~\ref{localreg} to the global case of the whole of $\Omega$.

\begin{proposition}[Global Hessian estimate]\label{stima hessiano globale}
Let $u\in C^{1}(\overline{\Om})$ be a solution to \eqref{eq debole}.
Then for $\beta\in [0,1)$ and $\gamma< n-2$ ($\gamma=0$ if $n=2$),
it holds
\begin{equation}\label{eq stima hessiano globale componente}
\sup\limits_{y\in \Omega}\, \io\frac{(k+|\n u|)^{p-2-\beta}|u_{ij}|^2}{|x-y|^\gamma}\di x \leq C
\end{equation}
and
\begin{equation}\label{eq stima hessiano globale}
\sup\limits_{y\in \Omega}\, \io\frac{(k+|\n u|)^{p-2-\beta}|D^2 u|^2}{|x-y|^\gamma}\di x \leq C\,,
\end{equation}
where $C= C(\beta,\gamma,p,n,\|u\|_{W^{1,\infty}},f)$ and $k\geq 0$ is given in Section \ref{intro}.
\end{proposition}
\proof
Recall that, since $\Om$ is smooth, then the boundary $\partial \Om$ satisfies the
interior sphere condition at each point. By Theorem \ref{hopf} we then have that $\n u$ does not vanish on $\partial \Om$.
Therefore by compactness of $\overline{\Om}$, applying Theorem \ref{stima hessiano locale}, we get the thesis.
\endproof

Arguing as in the proof of Proposition \ref{stima hessiano locale},
using once again Hopf's Lemma to ensure that $\n u$ does not vanish on the boundary of $\Omega$,
we can extend also the local result in Proposition \ref{inverso peso locale} to the following global result.

\begin{proposition}[Global estimate of the weight]\label{inverso peso globale}
Let $u\in C^{1}(\overline{\Om})$ be a solution to \eqref{eq debole}.
Fix $t\in [0,p-1)$ and $\gamma<n-2$ ($\gamma=0$ if $n=2$).
Then there exists $C>0$ such that
\begin{equation}\label{stima peso locale}
\sup\limits_{y\in \Omega}\,\, \int_{\Omega}\frac{\di x }{(k+|\n  u|)^t |x-y|^\gamma}\leq C\,,
\end{equation}
where $C=C(\Omega,t,\gamma,n,p,\|u\|_{W^{1,\infty}}, f)$. Moreover, $|Z| = 0$, where $Z$ is the critical set of $u$.
\end{proposition}

We finish this section with a discussion about Sobolev regularity for our positive solutions, as a corollary of the above global results.

\begin{theorem}\label{teo reg}
Let $u$ be a $C^{1}(\overline\Omega)$ solution
to $\eqref{eq debole}$.
There holds:
\begin{enumerate}
\item if $p\in (1,3)$, then $u\in W^{2,2}(\Om)$
\item if $p\in [3,+\infty)$, then $u\in W^{2,q}(\Om)$ with $q\leq\frac{p-1}{p-2}$.
\end{enumerate}
\end{theorem}
\proof
For any $p>1$, taking $\gamma=0$ in \eqref{eq stima hessiano globale} and \eqref{stima peso locale} respectively, we immediately get:
\begin{equation}\label{eq reg 01}
\io(k+|\n u|)^{p-2-\beta}|D^2 u|^2\di x \leq C\,,
\end{equation}
as well as
\begin{equation}\label{eq reg 02}
\int_{\Omega}\frac{1}{(k+|\n  u|)^t} \di x  \leq C\,,
\end{equation}
for any $t < p-1$.
Now, for $1<p<3$, choosing  $\beta < 1$ such that $p-2-\beta < 0$ and
recalling that $\n u$ is bounded, from \eqref{eq reg 01} we immediately get:
\begin{equation*}
\int_\Om |u_{ij}|^2 \di x  \leq \sup_\Om (k+|\n u|)^{\beta+2 - p} \int_\Om (k+|\n u|)^{p-2-\beta}  |u_{ij}|^2 \di x \leq  C,
\end{equation*}
which proves statement $(1)$. On the other hand, for $p \geq 3$, by \eqref{eq reg 01} and \eqref{eq reg 02} we have:
\begin{equation}\label{eq reg 1}
\begin{split}
&\int_\Om |u_{ij}|^q \di x =
  \int_\Om |u_{ij}|^q (k+|\n u|)^{(p-2-\beta)\frac{q}{2}}\cdot\frac{1}{(k+|\n u|)^{(p-2-\beta)\frac{q}{2}}} \di x\\
 &\leq
  \left(\int_\Om|u_{ij}|^2(k+|\n u|)^{p-2-\beta} \di x \right)^{\frac{q}{2}}
  \left(\int_\Om \frac{1}{(k+|\n u|)^{(p-2-\beta)\frac{q}{2-q}}} \di x\right)^{\frac{2-q}{2}} \leq C,
\end{split}
\end{equation}
where we used Holder's inequality with exponents
$\frac{2}{q}$ and $\frac{2}{2-q}$. Notice that, in order to apply \eqref{eq reg 02}, we need $(p-2-\beta) > 0$, which is true for $p \geq 3$, as well as
$(p-2-\beta)\frac{q}{2-q}<p-1$, and this is ensured taking
$q<\frac{p-1}{p-2}$ and recalling that $\beta<1$. This proves the statement $(2)$.

Finally, recalling that $|Z|=0$ and closely following the steps
in the proof of Proposition 2.2 in \cite{DS1}, we infer
that the generalized derivatives of $u_i$ coincide with the classical ones
almost everywhere in $\Om$ and we get the thesis.
\endproof

\begin{remark}
In case we cannot apply Hopf's Lemma, the same regularity
results stated in Theorem \ref{teo reg} can be locally obtained
using the local estimates \eqref{eq stima hessiano locale componente},
\eqref{eq stima hessiano locale} and \eqref{stima peso locale}.
\end{remark}

\end{document}